\documentclass[12pt]{article}

\usepackage{amsfonts}

\usepackage{color}

\newtheorem{thm}{Theorem}[section]
\newtheorem{lemma}[thm]{Lemma}
\newtheorem{cor}[thm]{Corollary}

\newtheorem{prop}[thm]{Proposition}

\newcommand{\proof
}{\par\medskip\noindent {\bf Proof.\ \ }}

\newcommand{\be}{\begin{equation}}
\newcommand{\ee}{\end{equation}}
\newcommand{\openbox}{\leavevmode
       \hbox to8pt{\hfil\vrule\vbox to6pt{\hrule width6pt\vfil\hrule}\vrule}}

\newcommand{\qed}{\hbox to5pt{ } \hfill \openbox\bigskip\medskip}

\newcommand{\Fq}{\mathbb F _q}
\newcommand{\ve}[1]{\mathbf{#1}}

\newcommand{\cT}{\mbox{$\cal T$}}

\newcommand{\cF}{\mbox{$\cal F$}}
\newcommand{\cG}{\mbox{$\cal G$}}
\newcommand{\cH}{\mbox{$\cal H$}}

\newcommand{\cM}{\mbox{$\cal M$}}

\newcommand{\N}{\mathbb N}

\newcommand{\Q}{\mathbb Q}

\newcommand{\F}{\mathbb F}

\title{Gr\"obner Bases for Increasing Sequences}

\author{G\'abor Heged\"{u}s\footnote{ \'Obuda University,
B\'ecsi \'ut 96, Budapest, Hungary,
H-1037,
{\tt hegedus.gabor@uni-obuda.hu}}
, Lajos R\'onyai\footnote{Institute for Computer Science and Control,
E\"otv\"os Lor\'and
Research Network; and Department of Algebra, Institute of Mathematics,
Budapest University of Technology and Economics, M\H {u}egyetem rkp. 3.,
H-1111 Budapest, Hungary,
{\tt lajos@info.ilab.sztaki.hu}}
}

\begin{document}

\footnotetext{The research of LR was supported in part by the Hungarian Ministry of
Innovation and Technology NRDI Office within the framework of the Artificial
Intelligence National Laboratory Program.}

\maketitle

\begin{abstract}
Let $q,n \geq 1$ be integers, $[q]=\{1,\ldots, q\}$, and
$\F$ be a field with $|\F|\geq q$.
The set
of increasing sequences
$$
I(n,q)=\{(f_1,f_2, \dots, f_n) \in [q]^n:~ f_1\leq f_2\leq\cdots \leq f_n \}
$$
can be mapped via an injective map $i: [q]\rightarrow \F $ into
a subset $J(n,q)$ of the affine space $\F^n$.
We describe reduced  Gr\"obner bases,
standard monomials and Hilbert function of the ideal  of polynomials
vanishing on $J(n,q)$.

As applications we give an interpolation basis for $J(n,q)$, and
lower bounds
for the size of increasing Kakeya sets, increasing Nikodym sets, and for
the size of affine hyperplane covers of $J(n,q)$.
\end{abstract}
\medskip

\noindent
\section{Introduction}

In the paper  $q,n \geq 1$ are integers, and $[q]=\{1,\ldots, q\}$. We view
$[q]$ as an ordered set with $1<2<\cdots <q$ and consider  the following set
of increasing sequences
$$
I(n,q)=\{(f_1,\ldots ,f_n) \in [q]^n:~~
f_1\leq f_2\leq\cdots \leq f_n \}.
$$
Clearly we have
$$
I(n,q)|={n+q-1 \choose q-1},
$$
as shown for example by a stars and bars argument.
Let $\F$ be a field with $|\F|\geq q$. We can map  $[q]$ into a  subset of
$\F$ with the aid of an injective map $i: [q]\rightarrow \F$. This induces
a map of
$I(n,q)$ to a subset $J(n,q)$ of the affine space $\F^n$: a
sequence $(v_1, \ldots , v_n)\in [q]^n$ is mapped to
$(i(v_1), \ldots , i(v_n))\in \F^n$.
One may study the algebraic and geometric properties of the set $J(n,q)$.
We refer to \cite{LHMO}, \cite{MR1}, \cite{MR2} and references therein
for research of this kind on other combinatorially relevant point sets and
polynomial ideals.

Most of the properties of $J(n,q)$ considered here are independent of the
map $i$, as long as it is an injection. One exception is grid maps. An
injection $i:[q]\rightarrow \F$ is a {\em grid map} if  there exists
an $a\in \F$ such that
$i(j)=a+j$ for $j\in [q]$. In particular, a grid map of $[q]$  exists iff
the characteristic of $\F$ is at least $q$.
Our main technical result is the determination of  Gr\"obner bases,
standard monomials and Hilbert functions (for the definitions please see
Section \ref{first}) of the ideal $\mathbf I(J(n,q))$ of polynomials
vanishing on $J(n,q)$
(Proposition
\ref{GB}, Corollary \ref{Hilbert} and a useful
interpolation basis for $J(n,q)$).
We obtained applications of these results in several directions, as it is
detailed next.

\subsection{Interpolation and covering by
hyperplanes}

We denote by $\F[x_1,\ldots, x_n]=\F[\ve x]$ the polynomial ring over $\F$
with variables $x_1,\ldots, x_n$.
The standard monomials for the ideal  $\mathbf I(J(n,q))$ form
a linear
basis of the space of functions form $J(n,q)$ to the ground field
$\F$. Here we exhibit an other natural basis of interpolation which turns
out to be useful when we consider coverings by hyperplanes. The polynomials
can be given quite explicitly when $i$ is a grid embedding.

\begin{thm} \label{main3}
Let $\ve s\in J(n,q)$. Then there exists a unique
polynomial $P_{\ve s}\in \F[\ve x]$ such that
\begin{itemize}

\item[(i)] $P_{\ve s}(\ve s)=1$ and $P_{\ve s}(\ve w)=0$ for
each $\ve w\in J(n,q)$, $\ve w\ne \ve s$, with $deg(P_{\ve s})=q-1$,
and
\item[(ii)] If $i$ is a grid embedding, then
we can write $P_{\ve s}$ into the form
$P_{\ve s}=\prod_{i=1}^{q-1} L_i$, where the $L_i$ are linear polynomials.
      \end{itemize}
\end{thm}

\medskip

There are several results on covering discrete subsets of $\F^n$ by
hyperplanes (see \cite{AF}, \cite{B}, \cite{BS}, \cite{BBS}, and \cite{Z}), a prominent
example being the theorem of Alon and F\"uredi  \cite{AF} on the covers of
discrete grids with the
exception of a point. An analogous, sharp statement for increasing vectors is the
following:

\begin{thm} \label{main44}
Let $0<k\leq n$ be an integer and let $\ve s_1,\ldots ,\ve s_k\in J(n,q)$ be
increasing vectors.
Let $\{H_j:~ 1\leq j\leq m\}$ be a set of affine hyperplanes  such that
$$
J(n,q)\setminus \{\ve s_1, \ldots ,\ve s_k\}\subseteq \cup_{j=1}^m H_j.
$$

Then $m\geq q-1$.
\end{thm}

A similar reasoning gives  a sharp (non-constant) bound on the number of
hyperplanes covering the whole $J(n,q)$.

\begin{thm} \label{main5}
Let $\{H_i:~ 1\leq i\leq m\}$ be a set of affine hyperplanes  such that
$$
J(n,q)\subseteq \cup_{i=1}^m H_i.
$$
Then
$m\geq q$.
\end{thm}

\subsection{Kakeya and Nikodym sets for increasing vectors}

Let  $\ve a, \ve v\in {\Fq}^n$ be vectors, with $\ve v\not= \ve 0$.
Define the
line $\ell(\ve a,\ve v)\subseteq \Fq ^n $ as
$$
\ell(\ve a,\ve v):=\{ \ve a + t\ve v:~ t\in  {\Fq}\}.
$$

A subset  $K\subseteq {\Fq}^n$ is a {\em Kakeya set}, if for
each $\ve 0 \not= \ve v\in {\Fq}^n$ there exists an $\ve a\in {\Fq}^n$ such
that $\ell(\ve a,\ve v)\subseteq K$.

\medskip

\medskip

Wolff proposed the conjecture in \cite{W}, that for every  $\epsilon >0$
and for every $n\geq 1$ there exists a constant $c(n,\epsilon)$ such that
for any Kakeya set $K\subseteq {\Fq}^n$ we have
$$
|K|\geq c(n,\epsilon) q^{n-\epsilon}.
$$

Dvir obtained a stronger bound in his breakthrough work \cite{D}:

\begin{thm} \label{dvir}
Let  $K\subseteq {\Fq}^n$ denote a Kakeya set. Then
$$
|K|\geq {q+n-1\choose n}.
$$
\end{thm}

In \cite{G} Ganesan obtained lower and upper bounds for the size of local
Kakeya sets, where the set of required directions of the lines is local in
the sense that it is possibly a proper subset of ${\Fq}^n\setminus \{0\}$.
We consider
here some special local Kakeya sets,
the so-called  increasing Kakeya sets.

A subset  $K\subseteq {\Fq}^n$ is an {\em increasing Kakeya set}, if
for each $\ve 0\not= \ve v\in J(n,q)$ there exists an $\ve a\in {\Fq}^n$ such
that the line $\ell(\ve a,\ve v)\subseteq K$.
Clearly each Kakeya set is an increasing Kakeya set as well.
We can prove the same lower bound as in Theorem \ref{dvir} for increasing
Kakeya sets. In fact, our results on the standard monomials and
the Hilbert function
of $J(n,q)$ allow
to use the argument of Dvir-Tao-Alon directly for increasing
Kakeya sets. This is a strengthening of Dvir's result, as we have a smaller
number of conditions (lines) to consider.

\begin{thm} \label{mon_kak}
Let  $K\subseteq {\Fq}^n$ be an increasing Kakeya set. Then
$$
|K|\geq {q+n-1\choose n}.
$$
\end{thm}

It was proved on page 3 of \cite{F}, that for each prime power $q$  there
exists a subset $K\subseteq \Fq ^2$ which is a union of $q$ lines with
different directions and with
$$
|K|={q+1\choose 2},
$$
providing an optimal increasing Kakeya set in the case $n=2$.
We also have an optimal construction for $n=q=3$, hence
the bound of Theorem  \ref{mon_kak} is sharp if $n=2$ or $q=2$ or $n=q=3$.

In general we have the following simple construction which seems to
be good for small values of $q$:

\begin{equation} \label{T-equation}
T(n,q):=\cup_{\ve 0 \not= \ve v\in J(n,q)} \ell(\ve 0, \ve v).
\end{equation}
Clearly we have
$$
|T(n,q)|\leq (q-1)\left(  {q+n-1\choose n}-(q-1)\right)+1.
$$

It would be interesting to see, if there is a general upper bound on the size
of the smallest increasing Kakeya sets, which is better than the
best available upper bound for general Kakeya sets.

\medskip

Nikodym sets are  closely related to Kakeya sets.
A subset $B\subseteq {\Fq}^n$ is a  {\em Nikodym set}, if for each
$\ve z\in {\Fq}^n $
there exists a line $\ell_{\ve z}\subseteq {\Fq}^n$ through $\ve z$
such that
$\ell_{\ve z}\setminus  \{\ve z\}\subseteq B$.
A variant of the Dvir-Alon-Tao-argument for Kakeya sets \cite{D}, (see
also Theorem 2.9 in \cite{Guth}) gives that the size of a Nikodym set
$B\subseteq \F_q^n$ is at least ${q+n-2 \choose n}$. Here we obtain a
version of this result for increasing vectors.
We shall consider a local kind of Nikodym sets.
$B\subseteq {\Fq}^n$ is an {\em increasing Nikodym set} if for
each $\ve z\in J(n,q)$
there exists a line $\ell_{\ve z}$ through $\ve z$ such that
$\ell_{\ve z}\setminus  \{\ve z\}\subseteq B$.

\begin{thm} \label{nico}
Let  $B\subseteq {\Fq}^n$ be an increasing Nikodym set. Then
$$
|B|\geq {n+q-2\choose n}.
$$
\end{thm}

The inequality of Theorem \ref{nico} strengthens the above bound for Nikodym sets
in that here we have fewer conditions (lines to take care of).
We note also that the set $T(n,q)$ from (\ref{T-equation}) is an increasing
Nikodym set.

\medskip

In Section 2 we collected the preliminaries about Gr\"obner  bases,
standard monomials and Hilbert functions. In Section \ref{GBproofs} we describe  Gr\"obner bases and  standard
monomials for the vanishing ideal of increasing
sequences $J(n,q)$.  We extend these
results to some special subsets of $J(n,q)$ and to strictly increasing
sequences. Applications are discussed in Section \ref{appli}.

\section{Notation and results from Gr\"obner theory}
\label{first}

A total ordering  $\prec$ on the monomials $x_1^{i_1}x_2^{i_2}\cdots
x_n^{i_n}$  composed from
variables $x_1,x_2,\ldots, x_n$ is a {\em term order}, if 1 is the
minimal element of $\prec$, and $\ve u \ve w\prec \ve v \ve w$ holds
for any monomials
$\ve u,\ve v,\ve w$ with $\ve u\prec \ve v$.
Important term orders are the lexicographic
order $\prec_l$ and the deglex order $\prec _{dl}$. We have
$$
x_1^{i_1}x_2^{i_2}\cdots x_n^{i_n}\prec_l x_1^{j_1}x_2^{j_2}\cdots
x_n^{j_n}
$$
iff $i_k<j_k$ holds for the smallest index $k$ such
that $i_k\not=j_k$. For deglex, we have $\ve u\prec_{dl} \ve v$ iff  either
$\deg \ve u <\deg \ve v$, or $\deg \ve u =\deg \ve v$, and $\ve u \prec_l
\ve v$.

The {\em leading monomial} ${\rm lm}(f)$
of a nonzero polynomial $f$ from the polynomial
ring $\F[\ve x]=\F[x_1, x_2, \ldots ,x_n]$ is the largest
monomial with respect to $\prec$,
which has nonzero coefficient in the standard form of $f$.

Let $I$ be an ideal of $ \F[\ve x] $. A finite subset $\cG\subseteq I$ is
a {\it
Gr\"obner basis} of $I$ if for every $f\in I$ there exists a $g\in G$ such
that ${\rm lm}(g)$ divides ${\rm lm}(f)$. It is known that such a $\cG$ is a basis of $I$.
A fundamental fact is (cf. \cite[Chapter 1, Corollary
3.12]{CCS} or \cite[Corollary 1.6.5, Theorem 1.9.1]{AL}) that every
nonzero ideal $I$ of $ \F[\ve x]$ has a Gr\"obner basis with respect to any term
order $\prec$.
\medskip

A monomial $\ve w\in \F[\ve x]$ is a {\em standard monomial} for $I$ if
it is not a leading monomial of any $f\in I$. Let ${\rm sm}(I,\prec)$ stand
for the
set of all standard monomials of $I$ with respect to the term-order
$\prec$ over $\F$.
It is known (see \cite[Chapter 1, Section 4]{CCS}) that for a
nonzero ideal $I$  the set
${\rm sm}(I,\prec)$ is a basis of the $\F$-vector space $ \F[\ve x]/I$.
In fact,
every $g\in \F[\ve x]$ can be written uniquely as $g=h+f$ where $f\in I$ and
$h$ is a unique $\F$-linear combination of
monomials from ${\rm sm}(I,\prec )$. For a subset $X\subseteq \F^n$ we write
$\mathbf I(X)$ for the ideal of polynomials from $ \F[\ve x]$ which vanish on $X$.
When $X\subseteq \F^n$ is a finite subset, interpolation
gives that every $X\rightarrow \F$ function
is a polynomial function. The
latter two facts imply for $\mbox{sm}(X,\prec):=\mbox{sm}(\mathbf I (X),\prec)$
that

\begin{equation}
\label{standard}
|\mbox{sm}(X,\prec)|=|X|.
\end{equation}

A Gr\"obner basis $\{ g_1,\ldots ,g_m\}$ of $I$ is {\em reduced}
if the coefficient of ${\rm lm}(g_i)$ is 1, and no nonzero monomial
in $g_i$ is divisible by any ${\rm lm}(g_j)$, $j\not= i$.
By a theorem of Buchberger (\cite[Theorem 1.8.7]{AL}) a nonzero ideal
has a unique reduced Gr\"obner basis.

The {\em initial ideal} ${\rm in}(I)$ of $I$ is the ideal in
$\F[\ve x]$ generated
by the monomials $\{ {\rm lm}(f):~f\in I\}$.

The notion of reduction is closely related to Gr\"obner bases.
Let $\cG$ be a Gr\"obner basis of and ideal $I$
of $\F[\ve x]$ and
     $f\in \F[\ve x ]$ be a polynomial.
We can reduce $f$ by the set $\cG$ by subtracting multiples of polynomials
$g\in \cG$ from $f$ in such a way that the resulting polynomial is composed
of $\prec$-smaller monomials.  It is known that this way any
$f\in \F[\ve x]$ can be
reduced into a (unique) $\F$-linear combination of standard monomials. This
is related to the fact we have already mentioned: ${\rm sm}(I,\prec )$
provides a linear basis of $\F[\ve x]/I$.

\medskip

Let $I$ be an ideal of $\F[\ve x]=\F[x_1,\ldots ,x_n]$. The (affine) {\em Hilbert
function} of the algebra $\F[\ve x]/I$ is the sequence of natural numbers
$h_{\F[\ve x]/I}(0), h_{\F[\ve x]/I}(1),
      \ldots $. Here $h_{ \F[\ve x]/I}(m)$ is the dimension over $\F$ of the
factor space
      $$\F[x_1,\ldots,x_n]_{\leq m}/(I\cap \F[x_1,\ldots,x_n]_{\leq m})$$
(see
\cite[Section 9.3]{BW}). It is easy to see that $h_{\F[\ve x]/I}(m)$ is the
number of standard monomials of degree at most $m$, where the ordering
$\prec$ is deglex.

For $I=\mathbf I (X)$ with some $X\subseteq \F ^n$, the number
$h_{X}(m):=h_{\F[\ve x]/I}(m)$ is the dimension of the linear space of those
$X\rightarrow \F$ functions which are polynomials of degree at
most $m$.

\section{Gr\"obner bases for increasing sequences}
\label{GBproofs}

Via an injective map $i: [q]\rightarrow \F$ we consider $[q]$ as a subset of
our ground field $\F$, in particular we assume that $|\F|\geq q$. We denote
by $J(n,q)$ the image of $I(n,q)$ by the map induced by $i$:
$$ J(n,q)=\{ (i(v_1), \ldots ,i(v_n)); ~~ (v_1,\ldots v_n)\in
I(n,q) \}.$$
The ordering $<$ on
$i([q])$ is defined via the map $i$, as $i(1)<i(2)< \cdots <i(q)$. Let $\prec$ be a term
order on the monomials of $\F[\ve x]$

We say that an $n$-tuple $(I_1,\ldots ,I_n)$ of subsets of $[q]$
is a good decomposition  of $[q]$, if

\begin{itemize}

\item[(i)] $\cup_{j=1}^n I_j=[q]$,
\item[(ii)] if $i<j$, then $x<y$ for each $x\in I_i,y\in I_j$.
\end{itemize}
In particular $I_k\cap I_l=\emptyset$
if $1\leq k\neq l\leq n$.

Next we define a polynomial $f_{I_1,\ldots ,I_n}(\ve x)\in \F[\ve x]$
attached to a good decomposition $(I_1,\ldots ,I_n)$ of $[q]$. We set
$$
f_{I_1,\ldots ,I_n}(x_1,\ldots x_n):= \prod_{j=1}^n
\left(\prod_{t\in i(I_j)} (x_j-t)\right)\in {\F}[x_1,\ldots x_n].
$$
Let
$$
\cG:= \{f_{I_1,\ldots ,I_n}:~ (I_1,\ldots ,I_n)
\mbox{ is a good decomposition of } [q] \}.
$$
\begin{prop} \label{GB}
$\cG$ is the reduced Gr\"obner basis of $J(n,q)$ for any term order $\prec$.
Moreover for the standard monomials we have
$$
{\rm sm}(J(n,q),\prec)= \{x^u:~ \mbox{deg}(x^u)\leq q-1 \}.
$$
\end{prop}

\proof
Observe first that the leading monomial of $f_{I_1,\ldots ,I_n}$ is
$x_1^{|I_1|}x_2^{|I_2|}\cdots x_n^{|I_n|}$ for any term order $\prec$ on
the monomials of $\F[\ve x]$, and
it is of degree $q$. Moreover any monomial of degree $q$ is the leading
monomial of a polynomial of the shape $f_{I_1,\ldots ,I_n}$ for a suitable
good decomposition $(I_1,\ldots ,I_n)$ of $[q]$.
It is then sufficient to prove that
$f_{I_1,\ldots ,I_n}(x_1,\ldots x_n)$ vanishes on the vectors of
$J(n,q)$.
Indeed then it follows that
$\mbox{sm}(J(n,q),\prec) \subseteq \{x^u:~ \mbox{deg}(x^u)\leq q-1\}$.
The two sets in the preceding formula have the same size
${n+q-1 \choose q-1}$, hence they must be equal:
$$\mbox{sm}(J(n,q),\prec) = \{x^u:~ \mbox{deg}(x^u)\leq q-1\}.$$
This gives the statement about the standard monomials and shows also
that the polynomials in $\cG$ indeed form a Gr\"obner
basis. These polynomials are monic, and except for the leading monomial
they are made of standard monomials. This proves that the Gr\"obner basis
$\cG$ is reduced.

It remains to verify that if $\ve v\in J(n,q)$, then $f_{I_1,\ldots
,I_n}(\ve v)=0$. Suppose that this is not the case,  $f_{I_1,\ldots
,I_n}(\ve v)\not=0$. Then straightforward induction on $j$ gives
that $v_j\not \in i(I_1\cup\cdots \cup I_j)$.
This leads in the end to $v_n\not\in i([q])$, a contradiction.
$\Box$

\medskip
Please note that the Gr\"obner basis and the standard monomials are
independent of the term order $\prec$ selected.
For the the Hilbert function of $J(n,q)$ we
have:

\begin{cor}\label{Hilbert}
$$
h_{J(n,q)}(s)={n+s \choose s}
$$
for each $0\leq s\leq q-1$.
\end{cor}

\proof We apply Proposition \ref{GB} with $\prec$ being the deglex order.
We obtain that the value of $h_{J(n,q)}(s)$ is the number of
monomials in $\F[\ve x]$ of degree at most $s$. $\Box$
\medskip

We obtain the following version of the Combinatorial Nullstellensatz \cite{A}
for increasing sequences.

\begin{cor} \label{mon_CombNsatz}
Let $0\not= f\in {\F}[x_1,\ldots x_n]$ be a polynomial
such that $deg(f)\leq q-1$. Then there
exists a vector $\ve v\in J(n,q)$ such that $f(\ve v)\ne 0$.
\end{cor}

\proof By Proposition \ref{GB} $f$ is a nontrivial linear combination of
standard monomials, hence it can not vanish on the entire set
$J(n,q)$. $\Box$.

\medskip

We now proceed to exhibit Gr\"obner bases for the ideals
$\mathbf I(\cal F)$
where $\emptyset\not= {\cal F}\subseteq J(n,q)$ is a
{\em downset} in the sense that if have $\ve u, \ve v\in J(n,q)$
with $\ve v\in {\cal F}$ and $\ve u\leq \ve v$ (component-wise) then
$\ve u \in {\cal F}$. We denote by $\cF ^c$ the complement of $\cF$
in $J(n,q)$.  We define first the
map
$$\phi:I(n,q)\to \{\ve v\in \{0,1,\ldots ,q-1\}^n:~
\sum_{i=1}^n v_i\leq q-1\}:$$
if $\ve f=(f_1,\ldots ,f_n)\in I(n,q)$ is an
increasing vector, then
$$
\phi(\ve f):= (f_1-1, f_2-f_1, \ldots ,f_n-f_{n-1}).
$$
It is straightforward to see that $\phi$ is a bijection.

Let $\ve g=(g_1,\ldots ,g_n)\in I(n,q)$ be an increasing
vector and consider
the following system of sets corresponding to $\ve g$:
$I_1(\ve g):=\{1,2,\ldots ,g_1-1\}$ and
$I_j(\ve g):=\{g_{j-1}, \ldots ,g_j-1 \}$
for each $2\leq j\leq n$. We define the following
polynomial
$$
f_{I_1(\ve g),\ldots ,I_n(\ve g)}(x_1,\ldots x_n):= \prod_{j=1}^n
\left(\prod_{t\in i(I_j(\ve g))}
(x_j-t)\right)\in {\F}[x_1,\ldots x_n].
$$
Here an empty product has value 1.
We have $\mbox{lm}(f_{I_1(\ve g),\ldots ,I_n(\ve g)})=\ve x^{\phi(\ve g)}$.
Set
$$
\cH:=\{(I_1(\ve g), \ldots, I_n(\ve g)):~ i(\ve g) \in {\cF}^c\},
$$
and
$$
\cT:=\{(I_1,\ldots ,I_n):~ (I_1,\ldots ,I_n)\mbox{ is a good decomposition
of } [q] \}.
$$

\begin{prop} \label{downset}
Let $\F$ be a field, $i([q])\subseteq \F$ as before, $\prec$ an arbitrary
term order on $\F[\ve x]$.
Let $\emptyset\not= \cF\subseteq J(n,q)$ be a downset. Then
$$
\cG:= \{f_{I_1,\ldots ,I_n}:~ (I_1,\ldots ,I_n)\in \cT \cup \cH \}
$$
is a Gr\"obner basis of the ideal $\mathbf I(\cF)$. Moreover for the
standard monomials we have
$$
{\rm sm}(\cF,\prec)=\{\ve x^{\ve u}:~ \ve u\in \phi(i^{-1}(\cF))\}.
$$
\end{prop}

Please note that the Gr\"obner basis $\cG$ is independent of the
term order $\prec$.

\proof
We observe first that the polynomials in $\cG$ all vanish on  $\cF$.
This was established for
$f_{I_1,\ldots ,I_n}$ when $(I_1,\ldots ,I_n)\in \cT$ in Proposition
\ref{GB}. Consider now vectors
$\ve v\in \cF$ and
$\ve g \in I(n,q)$ such that $i(\ve g)\in \cF^c$.
We have to establish $f_{I_1(\ve g),\ldots ,I_n(\ve g)}(\ve v)=0$.
If $f_{I_1(\ve g),\ldots ,I_n(\ve g)}(\ve v)\not =0$, then an induction
on $j$ gives that $v_j\geq i(g_j)$ holds for $j=1,\ldots ,n$. This is
immediate for
$j=1$, as $v_1$ can not be in $i(I_1(\ve g))$. Also for $j>1$ the facts
$i(g_{j-1})\leq v_{j-1}\leq v_j$ and $v_j\not\in i(I_j(\ve g))$ give
the claim.

On the other
hand, by the selection of $\ve v$ and $\ve g$ there must be an index $j$
such that $i(g_j)>v_j$, giving a contradiction, which shows that the
polynomials from $\cG$ indeed vanish on $\cF$; here we used that
$\cF$ is a downset.

Next we verify that any monomial $\ve x^{\ve u}$ such that
$u\not\in \phi(i^{-1}(\cF))$ is divisible by the leading monomial
of a polynomial
from $\cG$. Suppose first that $deg(\ve x^{\ve u})\geq q$. Then there
exists a good decomposition $(I_1,\ldots ,I_n)\in \cT$ and a
$\ve w \in {\N}^n$ such that
$\mbox{lm}(\ve x^{\ve w}\cdot f_{(I_1,\ldots ,I_n)})=\ve x^{\ve u}$.
Suppose now that $deg(\ve x^{\ve u})< q$. As $u\not\in \phi( i^{-1}(\cF))$,
we have  $\ve u\in \phi( i^{-1}({\cF}^c))$. Let $\ve g:=\phi^{-1}(\ve u)$ and
consider $f_{I_1(\ve g),\ldots ,I_n(\ve g)}$ whose leading monomial is
$\ve x^{\phi(\ve g)}=\ve x^{\ve u}$. We obtained
$$
\mbox{sm}(\cF,\prec) \subseteq \{\ve x^{\ve u}:~ \ve u\in \phi(i^{-1}(\cF))\}.
$$
We have equality here, as both sides have size $|\cF|$.
Also, with $\cG$
we can reduce any polynomial into a linear combination of standard
monomials. This implies that $\cG$ is a Gr\"obner basis. $\Box$

\medskip

We remark that the statement above implies also that if
$\cF\subseteq J(n,q)$ is a downset
then $\phi(i^{-1}(\cF))$ is a downset as well. This fact can be seen more
directly, without using Gr\"obner theory,  by the fact
that $\phi^{-1}$ is an order preserving map.

\medskip

Let $q\geq n\geq 1$ be  fixed integers. Next we consider the collection
of strictly increasing sequences. We put
$$
SI(n,q)=\{f\in [q]^n:~ f\mbox{ is strictly  increasing }\}.
$$
Clearly we have $|SI(n,q)|={q\choose n}$.

Similarly to increasing sequences,  we view $[q]$ as an ordered
subset of a field $\F$ via an injective map
$i:[q]\rightarrow \F$, in particular we assume $|\F|\geq q$.
We consider the image $SJ(n,q)$ of
$SI(n,q)$ with respect to $i$:
$$
SJ(n,q):=\{ (i(g_1),\ldots, i(g_n)): ~~ \ve g\in SI(n,q)\}. $$
The set $SJ(n,q)$ is a subset of the affine space $\F^n$,
and we proceed to determine a Gr\"obner basis for the ideal of
polynomials from $\F[\ve x]$ which vanish on $SJ(n,q)$.
Let $\prec$ denote a term order on the monomials of $\F[\ve x]$.

Let $1\leq j_1<j_2<\cdots < j_{n-1}\leq q$ be integers.
The collection
$(I_1,I_2,\ldots, I_n)$ of subsets of $[q]$ is called a {\em super
decomposition} of $[q]$ if $I_1=\{1,2, \ldots, j_1-1\}$,
$I_2=\{j_1+1,\ldots ,j_2-1\}$ and so on, finally
$I_n=\{j_{n-1}+1,\ldots, q \}$. If $j_{m+1}=j_m +1$, then $I_m=\emptyset$. Clearly the sets $I_i$ are mutually disjoint
and we have $|\cup _{i=1}^n I_i|=q-n+1$.

Let
$$
\cM:=\{(I_1,\ldots ,I_n):~ (I_1,\ldots ,I_n)\mbox{ is a
super decomposition of } [q] \}.
$$
For a super decomposition $(I_1,\ldots ,I_n)$ we consider the polynomial
$$
f_{I_1,\ldots ,I_n}(x_1,\ldots x_n):=
\prod_{j=1}^n  \left(\prod_{t\in i(I_j)} (x_j-t)\right)\in {\F}[x_1,\ldots
,x_n],
$$
and set
$$
\cG:= \{f_{I_1,\ldots ,I_n}:~ (I_1,\ldots ,I_n)\in \cM \}.
$$
We have a statement analogous to Proposition \ref{GB}, with a similar proof.

\begin{prop} \label{GB2}
$\cG$ is the reduced Gr\"obner basis of $SJ(n,q)$. Moreover
$$
{\rm sm}(SJ(n,q),\prec)=
\{\ve x^{\ve u}:~ \mbox{deg}(\ve x^{\ve u})\leq q-n \}.
$$
\end{prop}

\proof
We note first that for any monomial $\ve w\in \F[x_1,\ldots ,x_n]$
of degree $q-n+1$ there exists a
super decomposition $(I_1,\ldots, I_n)$ of $[q]$ such that the leading
monomial of $f_{I_1,\ldots ,I_n}$ is $\ve w$. We claim first that it
suffices to verify that the polynomials $f\in \cG$ all vanish on
$SJ(n,q)$. Indeed, then we obtain at once that
$$
\mbox{sm}(SJ(n,q),\prec)\subseteq
\{\ve x^{\ve u}:~ \mbox{deg}(\ve x^{\ve u})\leq q-n \},
$$
which implies that the sets on the two sides are equal because they have
the same size ${q \choose n}$. We conclude from here as in Proposition
\ref{GB}, and obtain that $\cG$ is a reduced Gr\"obner basis.

It remains to verify that for every
$(I_1,\ldots ,I_n)\in \cM $, and
$\ve v\in SJ(n,q)$ we have  $f_{I_1,\ldots ,I_n}(\ve v) =0$.
Suppose that $\ve v$ is a counterexample.
A straightforward induction on $\ell$ gives that if
$f_{I_1,\ldots ,I_n}(\ve v) \not=0$, then $v_\ell\geq i(j_\ell)$ for
$\ell =1,\ldots
,n-1$, where $j_1,\ldots ,j_{n-1}$ is the sequence defining the super
decomposition $I_1,\ldots ,I_n$. This implies that $v_n>v_{n-1}\geq
i(j_{n-1})$, and
hence $v_n\not\in i([q])$ a contradiction. This finishes the proof.
$\Box$

\medskip

Applying the preceding result to a deglex order provides the Hilbert function
of $SJ(n,q)$.

\begin{cor}\label{SIHilbert}
$$
h_{SJ(n,q)}(s)={n+s \choose s}
$$
for each $0\leq s\leq q-n$. $\Box$
\end{cor}

\medskip

Similarly to Corollary  \ref{mon_CombNsatz}
we have a non-vanishing statement here as well.

\begin{cor} \label{smon_CombNsatz}
Let $f\in {\F}[x_1,\ldots x_n]$ be a polynomial such
that $deg(f)\leq q-n$. Then there exists a $\ve v\in SJ(n,q)$ such that
$f(\ve v)\ne 0$. $\Box$
\end{cor}

\section{Applications}
\label{appli}

\subsection{Interpolation and covering}

\noindent
{\bf Proof of Theorem \ref{main3}.}
From Proposition \ref{GB} we obtain that the function
$P_{\ve s}: J(n,q)\rightarrow \F$ can be obtained as a unique
linear combination of standard monomials whose degree is at most
$q-1$. The degree of this polynomial $P_{\ve s}$ can not be smaller than
$q-1$, as otherwise the nonzero polynomial $(x_1-s_1)P_{\ve S}$ would
have degree $\leq q-1$ and vanish on $J(n,q)$, which is impossible by
Corollary \ref{mon_CombNsatz}. This proves (i).

We explain the proof for (ii) in the case when $\F=\Q$ and $i$
is the identical map of $[q]$. The general case follows similarly, but
involves more complicated notation. Let
$\ve s=(s_1,\ldots , s_n)\in J(n,q)\subseteq [q]^n$.
We shall give the linear factors
$L_i$ of $P_{\ve s}$ in terms of $\ve s$. We include the linear
polynomials $x_1-t$ for each $t\in [q]$ such that
$t<s_1$ and $x_n-t$ for each $t>s_n$. There
are as many as $q-(s_n-s_1+1)$ such linear factors. Moreover, for
each $i>1$ such that $s_i-s_{i-1}=k>0$, we consider the $k$
polynomials
$$
(x_i-x_{i-1}), (x_i-x_{i-1}-1), \ldots  ,(x_i-x_{i-1}-(k-1)).
$$
There are $s_n-s_1$ linear polynomials of this type. The product $Q$ of
the preceding $q-1$ linear polynomials does not vanish on $\ve s$. Suppose
that
$\ve t\in J(n,q)$ is a vector such that $Q(\ve t)\not=0$. Then an
inductive argument proceeding from $i=1$ to $i=n$ shows that
$\ve t\geq \ve s$ (component-wise comparison). A similar argument in
the direction from $i=n$ to $i=1$ gives $\ve t\leq \ve s$, hence
$\ve t=\ve s$. A suitable scalar multiple of $Q$ will be appropriate as
$P_{\ve s}$. Uniqueness follows as in the general case.
$\Box$

\medskip
\noindent

{\bf Example.} Let $n=5$,  $q=5$, and $\ve s=(1,2,2,4,4)$.
Then we have $Q(\ve x)=(x_5-5)(x_4-x_3)(x_4-x_3-1)(x_2-x_1)$.

\medskip

\noindent

{\bf Remarks.} 1. The polynomials $Q$ and $P_{\ve s}$ in (ii) are explicitly
determined by $\ve s$ in the proof. \\
      2. The uniqueness of $P_{\ve s}$ in (i) can also be established
by a dimension counting argument. The polynomials $P_{\ve s}$ form a basis the
space of $J(n,q)\rightarrow \F$ functions and hence also of
the space of polynomials of degree at most $q-1$ in $\F[\ve x]$. In
particular no nonzero polynomial from the latter space can be the
identically 0 function on $J(n,q)$. This reasoning gives also an
alternative proof of the part of Proposition  \ref{GB}
on standard monomials, when $\prec$ is the deglex order,
and then of Corollary \ref{Hilbert}.

\begin{lemma}   \label{upper_bound}
Let $0<k\leq n$ be an integer and let $\ve s_1,\ldots ,\ve s_k\in J(n,q)$
be increasing vectors.
Let $P\in \F[\ve x]$ be a polynomial such that $P(\ve s_i)\ne 0$
for each $i$ and $P(\ve w)=0$ whenever
$\ve w\in J(n,q)\setminus \{\ve s_1, \ldots ,\ve s_k\}$.
Then $deg(P)\geq q-1$.
\end{lemma}

\proof
Suppose for contradiction that there
exists a polynomial $P\in \F[\ve x]$ such that $P(\ve s_i)\ne 0$
for each $i$, but  $P(\ve w)=0$ for each
$\ve w\in J(n,q)\setminus \{\ve s_1, \ldots ,\ve s_k\}$,
and $deg(P)< q-1$. There exists a hyperplane in $\F^n$ which contains the
points $\ve s_i$, that is, a linear polynomial $L\in \F[\ve x]$ such
that $L(\ve s_i)=0$ for each $i$. We define the nonzero
polynomial $Q:=P\cdot L$.
Then $Q(\ve w)=0$ for each   $\ve w\in J(n,q)$. But
$\deg(Q)\leq q-1$, which contradicts to Corollary \ref{mon_CombNsatz}. \qed

\medskip

\noindent
{\bf Proof of Theorem \ref{main44}.}
Let  $L_j\in \F[\ve x]$ be a linear polynomial whose set of zeros is
exactly $H_j$.  We can apply Lemma \ref{upper_bound} to $P=L_1\cdots L_m$.
$\Box$

\medskip

\noindent
{\bf Proof of Theorem \ref{main5}.}
Let $L_i$ be a linear polynomial defining $H_i$. Then $P=L_1\cdots
L_m$ vanishes on $J(n,q)$ and $\deg P=m$. Here $m<q$ would contradict to Corollary \ref{mon_CombNsatz}.
$\Box$

\subsection{Increasing Kakeya and Nikodym sets}

The following statement is a variant of a result of Alon, Dvir and Tao,
see for example Theorem 1.5 in \cite{D}. We give the proof for the readers'
convenience. Here $\prec$ is an arbitrary term order on the variables
$x_1,\ldots ,x_n$.

\begin{prop} \label{main2} Let $q$ be a prime power,
Let $0< \ell \leq q-1$ be an integer.
Let $\cT\subseteq  \Fq ^n$ be a subset for which
\begin{equation} \label{star}
\{\ve x^{\ve b}:~ {\rm deg}(\ve x^{\ve b})\leq \ell \}
\subseteq {\rm sm}(\cT,\prec).
\end{equation}
Suppose that $K\subseteq \Fq ^n$ is a subset such that for each
$\ve 0\not= \ve v\in \cT$ there exists a vector  $\ve a\in {\Fq}^n$ such that
$|\ell(\ve a,\ve v)\cap K|\geq \ell +1$.
Then
$$
|K|\geq {n+\ell\choose n}.
$$
\end{prop}

\proof
Note first that by the assumption $\ell>0$ we have $|\cT|>1$, hence
$K$ is not empty. Suppose for contradiction that
$$
|K|<  {n+\ell\choose n}.
$$
Then the monomials from $\{\ve x^{\ve b}:~ {\rm deg}(\ve x^{\ve b})\leq \ell \}$
can not be linearly independent as functions on $K$:
there exists a nonzero polynomial $P\in \Fq[\ve x]$ such
that $D:={\rm deg}(P)\leq \ell$ and $P(\ve v)=0$ for each $\ve v\in K$. We
can write $P$ as a sum of two polynomials:
$$
P=P_D+Q,
$$
where $P_D\not=0 $ is the homogeneous part of $P$ of degree $D$ and
${\rm deg}(Q)<D$. Note that we have $D>0$, because otherwise $P$ would be a nonzero constant function. This implies that $P_D(\ve 0)=0$.

Let $\ve 0\not= \ve v\in \cT$. Then there exists a vector
$\ve a\in {\Fq}^n$
such that  we have $|\ell(\ve a,\ve v)\cap K|\geq \ell +1$.  Define now
a polynomial in the single variable $t$ as follows:
$$
P_{\ve v}(t):=P(\ve a+t \cdot \ve v)\in \Fq[t].
$$
The coefficient of $t^D$ in $P_{\ve v}$  is exactly $P_D(\ve v)$ and
${\rm deg}(P_{\ve v})\leq \ell$.

The condition $|\ell(\ve a,\ve v)\cap K|\geq \ell +1$ implies that there
exist different values
$t_1, \ldots , t_{\ell +1}\in \Fq$ such that  $P_{\ve v}(t_i)=0$ for
each $i$, implying that $P_{\ve v}$ is the identically zero polynomial.

From this we obtain that  $P_D(\ve v)=0$ for each  $\ve v\in \cT$.
This contradicts to (\ref{star}), as $P_D$ is  a nontrivial linear
combination of standard monomials, hence can not vanish on the
entire $\cT$. This finishes the proof. \qed

\medskip

\noindent
{\bf Proof of Theorem \ref{mon_kak}.} We apply Proposition \ref{main2}.
Set $\ell=q-1>0$ and $\cT=J(n,q)$. Proposition \ref{GB} shows that
condition (\ref{star}) is satisfied.  Let $K \subseteq {\F_q}^n$ be an
increasing Kakeya set. Proposition \ref{main2} applies and gives the lower
bound
$$ |K|\geq {n+q-1\choose n}. $$
$\Box$

\medskip

\noindent
{\bf Example.} We describe an increasing Kakeya set
$K\subset {\F_3}^3$ for the set $J(3,3)$
of nondecreasing vectors. Take the "usual" ordering $0<1<2$ of $\F_3$. The
10 nondecreasing vectors from ${\F_3}^3$ appear in two groups:

$\ve v^0=(0,0,0), ~~~ \ve v^1=(0,0,1), ~~~ \ve v^2=(0,0,2), ~~~
\ve v^3=(0,1,1)$, \\

$\ve v^4=(0,1,2), ~~~ \ve v^5=(0,2,2)$, \\
these are from the plane $L=\{x_1=0\}$, and there are four more points
outside $L$: \\

$\ve w^1=(1,1,1),  ~~~ \ve w^2=(1,1,2), ~~~ \ve w^3=(1,2,2),  ~~~
\ve w^4=(2,2,2)$.

Among these 10 vectors a maximal projectively inequivalent system
is $\ve v^1$,
$\ve v^3$, $\ve v^4$, $\ve w^1$, $\ve w^2$, $\ve w^3$. A Kakeya set for the
nondecreasing
vectors is then a union of six lines. In our construction we start out with
an optimal construction $K_0$ for the 3 lines in $L$. Note that this is
essentially the case of nondecreasing sequences in ${\F_3}^2$, hence
$|K_0|=6$. We select $K_0$ as the union of the three lines
$\{\lambda \ve v^1\}$,
      $\{\ve v^1+\lambda \ve v^3\}$, $\{ \lambda \ve v^4\}$. We have then

$$ K_0=\{(0,0,0),(0,0,1),(0,0,2),(0,1,2),(0,2,0),(0,2,1)\}.$$

We add to $K_0$ three more lines in directions $\ve w^1,\ve w^2,\ve w^3$
respectively,
which intersect the plane $L$ in $K_0$. These are \\

$\{(0,0,1)+\lambda \ve w^1\}$, \\

$\{(0,0,0)+\lambda \ve w^2\}$, \\

$\{(0,2,0)+\lambda \ve w^3\}$.

These three lines all pass through $(1,1,2)$; and together they include 3
      more points with first coordinate value 2. We conclude that the union $K$ of
      the six lines has $6+1+3=10$ points, hence $K$ is an optimal Kakeya set for
the nondecreasing vectors in ${\F_3}^3$.

\medskip

{\bf Proof of Theorem \ref{nico}.}
Suppose for contradiction that there exists an increasing Nikodym set $B$
with size
$$
|B|<{q+n-2 \choose n}.
$$
Then there exists  a nonzero polynomial $P\in \Fq[\ve x]$ such that
${\rm deg}(P)\leq q-2$ and $P(\ve v)=0$ for each $\ve v\in B$.

Let $\ve z\in J(n,q)\subset \F_q^n$ be an arbitrary element.
Then there exists a line $\ell_{\ve z}=\{\ve z+t\ve v:~ t\in \Fq\}$
with $\ve 0\not= \ve v\in \F_q^n$
through $\ve z$ such that $\ell_{\ve z}\setminus  \{\ve z\}\subseteq B$.

Define the polynomial $Q(t):=P(\ve z+t \ve v)\in \Fq[t]$.
Then $Q(t)=0$ for each $t\in ({\Fq})^*$,  because
$\ell_{\ve z}\setminus  \{\ve z\}\subseteq B$.
It follows from  ${\rm deg}(Q)\leq q-2$ that
$Q$ is the identically $0$  polynomial, hence  $P(\ve z)=Q(0)=0$.
We obtained that  $P(\ve z)=0$ for each $\ve z\in J(n,q)$,
and then
${\rm deg}(P)\leq q-2$ implies that
$P$ is the identically $0$  polynomial by Proposition \ref{GB}.
This contradiction proves the
claim. \qed



\begin{thebibliography}{MM}


\bibitem{A} N. Alon.  Combinatorial Nullstellensatz. {\em Comb.,
Probability and Computing}, {\bf 8(1-2)}  (1999), 7-29.

\bibitem{AF} N. Alon, Z. F\"uredi. Covering the cube by
affine hyperplanes. {\em European Journal of Combinatorics}, {\bf 14(2)}
(1993), 79-83.

\bibitem{AL} W. W. Adams, P. Loustaunau. {\it An Introduction to
Gr\"obner Bases,} American Mathematical Society, 1994.

\bibitem{B} S. Ball. On intersection sets in Desarguesian affine spaces.
{\em European Journal of Combinatorics}, {\bf 21(4)} (2000), 441-446.

\bibitem{BW} T. Becker, V. Weispfenning. {\em Gr\"obner bases - a
computational approach to commutative algebra,} Springer-Verlag, Berlin,
Heidelberg, 1993.

\bibitem{BBS} A. Blokhuis, A.E. Brouwer, T. Sz\H {o}nyi. Covering all
points except one.
{\em Journal of Algebraic Combinatorics,} {\bf 32(1)} (2010), 59-66.

\bibitem{BS} A.E. Brouwer, A. Schrijver. The blocking number of an affine space.
{\em Journal of Combinatorial Theory, Series A,} {\bf 24(2)} (1978), 251-253.

\bibitem{CCS} A. M. Cohen, H. Cuypers, H. Sterk (eds.). {\it Some
Tapas of Computer Algebra,} Springer-Verlag, Berlin, Heidelberg, 1999.

\bibitem{CLS} D.~Cox, J.~Little and D.~O'Shea.
{\it Ideals, Varieties, and Algorithms,} Springer-Verlag, Berlin,
Heidelberg, 1992.

\bibitem{D} Z. Dvir. On the size of Kakeya
sets in finite fields. {\em Journal of the American Math.
Soc.,} {\bf 22} (2009), 1093-1097.

\bibitem{DKSS} Z. Dvir, S. Kopparty, S. Saraf, M. Sudan. Extensions
to the method of
multiplicities, with applications to Kakeya sets and mergers.
{\em SIAM J. Comput.,}  {\bf 42(6)} (2013), 2305-2328.

\bibitem{F} X.W.C. Faber. On the finite field Kakeya problem
in two dimensions. {\em Journal of Number Theory,}
{\bf 124(1)} (2007), 248-257.


\bibitem{G} G. Ganesan.  Size of local finite field Kakeya sets.
{\em In: Extended Abstracts EuroComb 2021,} (pp. 1-4).
Birkh\"auser, Cham.

\bibitem{GKS} A. Guo, S.  Kopparty, M. Sudan. New affine-invariant
codes from lifting. {\em  Proceedings of the 4th conference on
Innovations in Theoretical Computer Science} (pp. 529-540) (2013).

\bibitem{Guth} L. Guth. {\it Polynomial methods in combinatorics,}
American Mathematical Soc., 2016.



\bibitem{LHMO} J.A. De Loera, C.J. Hillar, P.N. Malkin, M. Omar.
Recognizing graph theoretic properties with polynomial ideals.
{\em Electronic Journal of Combinatorics,} {\bf 17} R114 (2010).



\bibitem{LSW} B. Lund, S. Saraf, C. Wolf.  Finite field Kakeya and
Nikodym sets in three dimensions. {\em SIAM Journal on Disc.
Math.,} {\bf  32(4)} (2018), 2836-2849.


\bibitem{MR1}
      T. M\'esz\'aros, L. R\'onyai.
Some combinatorial applications of Gr\"obner bases. {\em In: Winkler, F. (eds)
Algebraic Informatics. CAI 2011.} Springer LNCS 6742,  pp. 65-83. (2011)


\bibitem{MR2} T. M\'esz\'aros, L. R\'onyai.
Standard monomials and extremal vector systems.
{\em Electronic Notes in Discrete Math.,} {\bf 61} (2017),  855-861.


\bibitem{SS}  S. Saraf, M. Sudan. An improved lower bound on
the size of Kakeya sets over
finite fields. {\em Analysis and  PDE,} {\bf 1(3)} (2008), 375-379.

\bibitem{W} T. Wolff. Recent work connected with the Kakeya problem.
{\em Prospects in Mathematics (Princeton, NJ, 1996),} {\bf 2} (1999),
129-162.

\bibitem{Z} C. Zanella. Intersection sets in $AG(n,q)$ and a
characterization of the
hyperbolic quadric in $PG(3,q)$. {\em Discrete Mathematics,} {\bf 255(1-3)}
(2002), 381-386.

\end{thebibliography}
\end{document}